\documentclass[11pt,oneside]{amsart}

\title{Infinite time decidable equivalence relation theory}
\author{Samuel Coskey}
\address{Samuel Coskey, The Graduate Center of The City University of
  New York, Mathematics Program, 365 Fifth Avenue, New York, NY 10016}
\email{scoskey@math.rutgers.edu, http://math.rutgers.edu/$\sim$scoskey}
\author{Joel David Hamkins}
\address{Joel David Hamkins, The Graduate Center of The City
  University of New York, Mathematics Program, 365 Fifth Avenue, New
  York, NY 10016 \& College of Staten Island of CUNY, Mathematics,
  2800 Victory Boulevard, Staten Island, NY 10314}
\email{jhamkins@gc.cuny.edu, http://jdh.hamkins.org}

\usepackage[marginratio=1:1]{geometry}
\usepackage{mydefs}
\usepackage{mathpazo,bm}
\usepackage{setspace}
\usepackage{tikz}
\usepackage{diagrams}
\newarrow{Blank}{}{}{}{}{}\def\SD{\rdBlank~{\subset}}\def\SU{\ruBlank~{\subset}}
\usepackage{graphicx}

\newcommand{\DD}{{\bm{D}}}
\newcommand{\DsD}{\bm{sD}}
\newcommand{\DE}{\bm{E}}
\newcommand{\DsE}{\bm{sE}}
\newcommand{\Coll}{\mathop{\mathrm{Coll}}}
\newcommand{\Abs}{\mathop{\mathrm{Abs}}}

\begin{document}
\begin{abstract}
  We introduce an analog of the theory of Borel equivalence relations
  in which we study equivalence relations that are decidable by an
  infinite time Turing machine.  The Borel reductions are replaced by
  the more general class of infinite time computable functions.  Many
  basic aspects of the classical theory remain intact, with the added
  bonus that it becomes sensible to study some special equivalence
  relations whose complexity is beyond Borel or even analytic.  We
  also introduce an infinite time generalization of the countable
  Borel equivalence relations, a key subclass of the Borel equivalence
  relations, and again show that several key properties carry over to
  the larger class.  Lastly, we collect together several results from
  the literature regarding Borel reducibility which apply also to
  absolutely $\dD2$ reductions, and hence to the infinite time
  computable reductions.
\end{abstract}
\maketitle

\section{Introduction}

The subject of Borel equivalence relation theory---by now a highly
developed, successful enterprise---begins with the observation (see
Friedman-Stanley \cite{friedman} and Hjorth-Kechris
\cite{hjorthkechris}) that many classification problems arising
naturally in mathematics can be regarded as relations, often Borel
relations, on a standard Borel space. For example, since groups are
determined by their multiplication functions, the isomorphism relation
on countable groups can be regarded as a relation on the subspace of
$2^{\omega\times\omega\times\omega}$ corresponding to the graphs of
such functions.  This isomorphism relation is properly analytic, but
its restriction to finitely generated groups is Borel. The subject
aims to understand these and many other relations by placing them into
a hierarchy of relative complexity measured by Borel
reducibility. Specifically, an equivalence relation $E$ on a Borel
space is \emph{Borel reducible} to another, $F$, if there is a Borel
function $f$ such that
\[x\mathrel{E}y\iff f(x)\mathrel{F}f(y)
\]
for all $x,y$ in the underlying space. In this case, we write
$E\leq_BF$, and we think of this reducibility as asserting that $F$ at
least as complex as $E$.  Indeed, the function $[x]_E\mapsto [f(x)]_F$
is an explicit classification of the $E$-equivalence classes using
$F$-equivalence classes.  More generally, composition with $f$ gives
an explicit method of turning any $F$-invariant classification into an
$E$-invariant classification, and in this sense, the classification
problem for $F$ is at least as hard as the classification problem for
$E$.

In this article, a small project, we aim to extend the analysis from
the Borel context to a larger context of effectivity.  Namely, we
shall consider the context of infinite time computability, a realm
properly between Borel and $\dD2$.  Specifically, we shall enlarge the
reducibility concept by allowing infinite time computable reduction
functions (a class of functions we review in
Section~\ref{Section.ITTMs}).  This is sensible for several
reasons. First, the class of infinite time computable functions
properly extends the class of Borel functions---a function is Borel
exactly when it is infinite time computable in uniformly bounded
countable ordinal time---while retaining much of the effective flavor
and content of the Borel context. The infinite time computable
functions, determined by the operation of a finite Turing machine
program computing in transfinite ordinal time, seem in many ways as
``explicit'' as the Borel functions are sometimes described to be, but
they reach much higher into the descriptive set-theoretic hierarchy,
well into the class $\dD2$.  Second, meanwhile, many natural
equivalence relations that lay outside the Borel boundary,
particularly those having to do with well-orders or with more
arbitrary isomorphism relations for countable structures, are captured
within the bounds of infinite time computability. For example, it is
infinite time computable, but not Borel, to decide whether a given
real codes a well order, and the corresponding order isomorphism
relation on countable well orders is infinite time computable, but not
Borel. More generally, the isomorphism relation for arbitrary
countable structures in arbitrary countable languages is also infinite
time computable, but not Borel. Third, it will turn out that much (but
not all) of the Borel theory carries over to our enlarged context, at
least for many of the relations studied by that theory.  Positive
instances of Borel reducibility, of course, carry over directly
because Borel functions are infinite time computable. Conversely, a
deep aspect of the Borel theory is that many of the proofs of
non-reducibility, that is, instances of equivalence relations $E$ and
$F$ for which $E\not\leq_B F$, actually overshoot the Borel context by
showing, for example, that there are no measurable reduction functions
for a given pair of equivalence relations; since infinite time
computable functions are measurable, these arguments also rule out
reducibilities in our context. The point is that such non-reduction
arguments lay at the center of the Borel theory, and the overshooting
phenomenon means that in these instances, the non-reduction results
transfer largely intact to the infinite time computable context.
Thus, in our project we explore the limits of this phenomenon. In
summary, we propose to study the hierarchy of equivalence relations
under the concepts of reducibility provided by infinite time
computability.

In contrast, recent work of Knight \cite{knight} and others aim far in
the other direction, by restricting the notion of reducibility from
Borel functions down to the class of (ordinary) Turing computable
functions.  Since this is a very restrictive notion of reduction, it
allowed them to separate many complexity classes more finely and to
analyze even classes of finite structures.  Our work here can be seen
as complementary to theirs, since we consider comparatively generous
notions of reducibility.

This paper is organized as follows.  In the second section, we shall
describe in detail the infinite time Turing machines, the class of
sets which they decide, and the class of functions which they
compute. In fact, we shall define several distinct ways in which these
machines can accept their input, leading to distinct but closely
related notions of effectivity. In the third section, we give some
basic facts about Borel equivalence relations and Borel reductions,
and compare this situation with the case of infinite time computable
equivalence relations and functions. In the fourth section, we
consider the special case of \emph{countable} Borel equivalence
relations (\emph{i.e.}, those with countable classes). We define the
class of infinite time enumerable equivalence relations, which is a
natural generalization of the class of countable Borel equivalence
relations to the infinite time context. In the last section, we give a
survey of methods of demonstrating \emph{non}-reducibility results in
the case of absolutely $\dD2$ reductions.  We use these methods to
determine the relationships between the infinite time computable
equivalence relations discussed in this paper, and these relationships
are summarized in Figure \ref{fig_diagram} at the conclusion of the
paper.

\section{The infinite time complexity classes}\label{Section.ITTMs}

Infinite time Turing machines, introduced in \cite{machines},
generalize the operation of ordinary Turing machines into transfinite
ordinal time. An infinite time Turing machine has three one-way
infinite tapes (the input tape, the work tape and the output tape),
each with $\omega$ many cells exhibiting either $0$ or $1$, and
computes according to a finite program with finitely many
states. Successor stages of computation proceed in exactly the
classical manner, with the machine reading from and writing to the
tape and moving the head left or right according to the program
instructions for the current state. At limit time stages, the
configuration of the machine is determined by placing the head on the
left-most cell, setting the state to a special ``Limit'' state, and
updating each cell of the tape with the $\limsup$ of the previous
values exhibited by that cell (which is the limit value, if the value
had stabilized, otherwise $1$). Computation stops only when the
``Halt'' state is explicitly attained, and in this case, the machines
outputs the contents of the output tape.  Since the input and output
tapes naturally accommodate infinite binary sequences, that is,
elements of Cantor space $2^\omega$, the machines provide infinitary
notions of computability and decidability on Cantor space. The
machines can be augmented with additional input tapes to accommodate
real parameters or oracles. We denote by $\varphi_e^z(x)$ the output
of program $e$ on input $x$ with parameter $z$, if this computation
halts; if it doesn't halt, then $\varphi_e^z(x)$ is undefined, or
diverges. A partial function $f\partialfrom
2^\omega\rightarrow2^\omega$ is infinite time \emph{computable} if
there exists $e$ and $z$ such that $f=\varphi_e^z$. For a program $e$
operating on a machine without a parameter tape, we denote the output
as $\varphi_e(x)$, and say that the corresponding function $\varphi_e$
is infinite time computable \emph{without parameters}.

A simple cofinality argument (\cite[Theorem~1.1]{machines}) shows that
if an infinite time computation halts, then it does so in a countable
ordinal number of steps. And if a computation does not halt, then it
is necessarily due to the fact that it is caught in an infinite loop,
in the strong sense that at limits of repetitions of this loop, the
computation remains inside the loop. (Looping phenomenon in ordinal
time is complicated by the curious fact that an infinite time
computation can exactly repeat a configuration, looping $\omega$ many
times, but nevertheless escape the loop at the limit of these
repetitions.)

A subset $A\subset2^\omega$ is infinite time \emph{decidable} if its
characteristic function is infinite time computable, and infinite time
\emph{semidecidable} if it is the domain of an infinite time
computable function.  These concepts naturally extend to subsets of
$(2^\omega)^n$ for $n\leq\omega$, by the use of canonical pairing
functions (or by using extra input tapes). We warn the reader that a
function can have an infinite time decidable graph, as a subset of the
plane, without being an infinite time computable function (see
\cite[Lost Melody theorem~4.9]{machines}). The reason has to do with
the inability of the infinite time machines to search entirely through
Cantor space in a computable manner, and so the analogue of the
classical algorithm showing this doesn't happen for finite-time
computations on $\omega$ simply does not work here. We say that a
function $f$ is infinite time \emph{semicomputable} if the graph of
$f$ is infinite time decidable. Thus, every every infinite time
computable function is infinite time semicomputable, but not generally
conversely.

We let $\DD$ denote the class of infinite time decidable subsets of
$2^\omega$. Since we have allowed a real parameter $z$ in the
definition of an infinite time decidable set, the class $\DD$ fits
naturally into the bold-face descriptive set theory context.
Similarly, let $\DsD$ denote the class of semidecidable subsets of
$2^\omega$, and $\widetilde{\DsD}$ the class of co-semidecidable
subsets of $2^\omega$.  The classes of infinite time decidable sets
and functions subsume the corresponding Borel classes.  In fact, we
have the following remarkable characterization.

\begin{thm}
  \label{thm_Borel_char}
  A set $A$ is Borel if and only if it is decided by a program which,
  for some countable ordinal $\alpha$, halts uniformly on all inputs
  in fewer than $\alpha$ steps.  Similarly, a function $f$ is Borel if
  and only if it is computed by a program which, for some countable
  ordinal $\alpha$, halts uniformly in fewer than $\alpha$ steps.
\end{thm}

Thus, Borel effectivity is the entry into the larger hierarchy of
computability provided by infinite time Turing machines.  Before the
proof, let us fix the notation $\omega_{1,ck}$ for the supremum of the
ordinals that are computable by a classical Turing machine.  If
$z\in2^\omega$, then let $\omega_{1,ck}^z$ denote the relativized
version, the supremum of the ordinals that are computable by a
classical Turing machine with oracle parameter $z$. The ordinal
$\omega_{1,ck}^z$ is known to be the least ordinal which is admissible
in the parameter $z$.  The proof of Theorem~\ref{thm_Borel_char} is
then an easy consequence of \cite[Theorem 2.7]{machines}, which states
that $A$ is $\Delta_1^1$ if and only if it can be decided by a program
which halts uniformly in fewer than $\omega_{1,ck}$ steps.

\begin{proof}[Proof of Theorem \ref{thm_Borel_char}]
  To establish the first statement, we need only observe that the
  proof of \cite[Theorem 2.7]{machines} relativizes to show that a set
  $A$ is $\Delta_1^1(z)$ if and only if it can be decided by a program
  which halts uniformly in fewer than $\omega_{1,ck}^z$ steps. This
  implies directly that every Borel set is infinite time decidable in
  uniformly bounded time. Conversely, if $A$ is infinite time
  decidable from parameter $z$ in time uniformly bounded by $\alpha$,
  then simply augment $z$ with a real $y$ coding $\alpha$, so that
  $\alpha<\omega_{1,ck}^{z\oplus y}$, and conclude that
  $A\in\Delta^1_1(z\oplus y)$ and therefore $A$ is Borel.

  For the second statement, suppose that $f$ is infinite time
  computable uniformly in some number of steps which is bounded below
  $\omega_1$. It follows that the graph of $f$ is infinite time
  decidable in some bounded number of steps, and therefore, by the
  first paragraph, the graph of $f$ is Borel. Conversely, if $f$ is a
  Borel function, then for each $n$, the set
  \[A_n\defeq\{x\in2^\omega:f(x)(n)=0\}
  \]
  is Borel, and hence it is infinite time decidable by a program which
  halts uniformly in fewer than $\alpha_n$ steps. It follows that
  $f(x)$ is infinite time decidable by a program which halts uniformly
  in fewer than $\sup(\alpha_n)+\omega$ steps.
\end{proof}

In the following proposition, $\Abs\dD2$ denotes the class of
absolutely $\dD2$ sets, where a set $A$ is \emph{absolutely} $\dD2$
when it is defined by a $\dP2$ formula $\phi$ and by a $\dS2$ formula
$\psi$, such that the formulas $\phi,\psi$ remain equivalent in any
forcing extension.

\begin{prop}
  \label{prop_containments}
  The classes of infinite time decidable, semidecidable, and
  co-semidecidable sets lie within the projective hierarchy as
  follows.
  \begin{displaymath}
  \begin{diagram}[h=.5em,w=2em]
    \dS1 &     &     &     & \DsD  &     &               \\
    {}   & \SD &     & \SU &       & \SD &               \\
    {}   &     & \DD &     &       &     & \;\;\Abs\dD2 \subseteq\dD2     \\
    {}   & \SU &     & \SD &       & \SU &               \\
    \dP1 &     &     &     & \widetilde{\DsD}
                                   &     &               \\
  \end{diagram}
  \end{displaymath}
\end{prop}

\begin{proof}
  That $\dP1$ sets are infinite time decidable follows from the fact
  that an infinite time Turing machine can detect whether a given
  relation is wellfounded (see \cite[Count-Through
  Theorem]{machines}).  That every $\DsD$ set is $\dD2$ is shown in
  \cite[Complexity Theorem]{machines}, but we briefly sketch the
  argument.  The idea is that any run of an infinite time computation
  can be coded by a real, namely, a code for a well-ordered sequence
  $\seq{r_\alpha}$, where each $r_\alpha$ is just a code for the
  configuration of the machine at stage $\alpha$.  It is $\dP1$ to
  check that a given real codes a well-order, and hence it is $\dP1$
  to check that a given real codes a computation history.

  Now, if $A$ is semidecidable, then for some program $e$ we have that
  $x\in A$ if and only if there exists a real code for a halting
  computation history for $e$ on input $x$, and hence $A\in\dS2$.  But
  also $x\in S$ if and only if \emph{every} real coding a
  \emph{settled} run of the program $e$ on input $x$ shows that it
  accepts $x$.  (A snapshot is said to be settled if and only if it
  shows that the program halts or is caught in a strongly repeating
  infinite loop and hence cannot halt.)  This shows that $A\in\dP2$,
  and since our $\dS2$ and $\dP2$ descriptions of $A$ are absolutely
  equivalent, we have that $\DsD\subset\Abs\dD2$.
\end{proof}

Note that the inclusions in Proposition~\ref{prop_containments} are
proper (except, consistently, the last), since by the classical
diagonalization argument the halting set
\[H\defeq\{(x,p):p\textrm{ halts on input }x\}
\]
is infinite time semidecidable but not infinite time decidable.  It
follows that the complement $\RR\times\omega\smallsetminus H$ is a set
which is absolutely $\dD2$ but not semidecidable.  It follows that the
function which maps each $x\in2^\omega$ to its infinite time
\emph{(light-face) jump}
\[x^\triangledown\defeq\{p:p\textrm{ halts on input x}\}
\]
is not infinite time computable.  Although one might expect that the
jump function is semicomputable, we shall see shortly that this is not
the case.  However, there is a program which \emph{eventually} writes
the jump, in the sense that on input $x$ the program will write
$x^\triangledown$ on the output tape and never change it after some
ordinal stage.  Indeed, consider the ``universal'' program which
simulates all programs simultaneously on the input $x$.  Each time one
of the simulated programs halts, the master program adds a code for
that program to a list on the output tape.  Since each halting program
will do so in countably many steps, the output tape will eventually
converge to a code for $x^\triangledown$.  We are thus led to study
the following broader classes of infinite time effective sets and
functions.

\begin{defn}\ 
  \begin{itemize}
  \item A partial function $f\partialfrom2^\omega\into2^\omega$ is
    infinite time \emph{eventually computable} if there exists a
    program $e$ such that on any input $x\in\dom(f)$, the computation
    of $e$ on $x$ has the feature that from some ordinal time onward,
    the output tape exhibits the value $f(x)$, and for
    $x\notin\dom(f)$, the output tape does not eventually stabilize in
    this way.
  \item A subset $A\subset 2^\omega$ is infinite time \emph{eventually
      decidable} if its characteristic function is infinite time
    eventually computable. We let $\DE$ denote the class of infinite
    time eventually decidable sets.
  \item A subset $A\subset 2^\omega$ is infinite time
    \emph{semieventually decidable} if it is the domain of an infinite
    time eventually computable function.  Denote by $\DsE$ the class
    of semieventually decidable sets and by $\widetilde{\DsE}$ the
    class of infinite time co-semieventually decidable sets.
  \end{itemize}
\end{defn}

Unlike the semicomputable functions, it is easy to see that the class
of infinite time eventually computable functions is indeed closed
under composition.  The class of infinite time eventually computable
functions retains many of the descriptive properties of the infinite
time computable functions, and as we have hinted, it contains some
useful non-infinite time computable functions.

\begin{prop}
  \label{prop_morecontainments}
  We can now extend Proposition~\ref{prop_containments} to show the
  containments among these new classes of subsets of $2^\omega$.
  \begin{displaymath}
  \begin{diagram}[h=.5em,w=1.5em]
    \dS1 &     &     &     & \DsD  &     &     &     & \DsE  &     &  \\
         & \SD &     & \SU &       & \SD &     & \SU &       & \SD &  \\
         &     & \DD &     &       &     & \DE &     &       &     &\;\;\Abs\dD2\subseteq\dD2\\
         & \SU &     & \SD &       & \SU &     & \SD &       & \SU &  \\
    \dP1 &     &     &     & \widetilde{\DsD}
                                   &     &     &     & \widetilde{\DsE}
                                                             &     &  \\
  \end{diagram}
  \end{displaymath}
  Each of these containments (except, consistently, the last) is
  proper.  Moreover, we have that $\DsE\cap\widetilde{\DsE}=\DE$.
\end{prop}

\begin{proof}
  Suppose that $A$ is semidecidable and let $e$ be a program which
  halts if and only if $x\in A$.  Let $q$ be the program which
  initially writes $0$, and then simulates $e$, changing its output to
  $1$ whenever $e$ halts.  Then $q$ converges to $1$ if $e$ halts and
  to $0$ if $e$ does not, and so $A$ is infinite time eventually
  decidable.  To see that every infinite time semieventually decidable
  set is absolutely $\dD2$, use the same argument as Proposition
  \ref{prop_containments}, but replace the halting notion of
  acceptance with eventual convergence, which is observable in the
  settled snapshot sequences.

  To see that the inclusions are proper, consider the following analog
  of the halting set.  Namely, let $S$ denote the ``stabilizing'' set
  $\{(x,p):p\textrm{ stabilizes on input }x\}$.  Then $S$ is easily
  seen to be infinite time semieventually decidable but not infinite
  time eventually decidable.  It follows that $S^c$ is absolutely
  $\dD2$ but not infinite time semieventually decidable.

  Finally, suppose that both $A$ and $A^c$ are $\DsE$.  Let $e$ be a
  program which eventually stabilizes if and only if the input $x\in
  A$ and let $q$ be a program which eventually stabilizes if and only
  if $x\not\in A$.  Then consider the program $r$ which simulates both
  $e$ and $q$, writing $1$ whenever $q$ changes its output, and
  writing $0$ whenever $e$ changes its output.  $r$ does not change
  its output until either of these events occurs.  Clearly, $r$ will
  eventually write $1$ if and only if $x\in A$ and it will eventually
  write $0$ if and only if $x\not\in A$.
\end{proof}

The relationship between the the corresponding classes of functions is
slightly different.

\begin{prop}
  A function $f$ is infinite time computable if and only if it is both
  infinite time eventually computable and infinite time
  semicomputable.
\end{prop}

It follows that the jump function $x\mapsto x^\triangledown$ is indeed
not semicomputable.

\begin{proof}
  If $f$ is infinite time eventually computable by program $e$ and
  semicomputable by program $q$, then it can be computed by the
  program which simulates $e$, at each step using $q$ to check to see
  if the value on the output tape for $e$ is correct.
\end{proof}

We have already observed that even the infinite time semieventually
decidable sets lie within the class of absolutely $\dD2$ sets.  When
it comes to functions, the absolutely $\dD2$ property only extends to
the infinite time eventually computable functions.  Here, a function
$f$ is said to be absolutely $\dD2$ if and only if its \emph{diagram}
\[\diag(f)\defeq\set{(x,s)\in2^\omega\times2^{\ord<\omega}\mid s\subset f(x)}
\]
is absolutely $\dD2$.  Not every function with an infinite time
decidable graph will be absolutely $\dD2$ in this sense, and indeed by
Theorem~\ref{thm_semi} in the next section, not every semicomputable
function is absolutely $\dD2$.

\begin{thm}
  \label{thm_provd12}
  Every infinite time eventually computable function is absolutely
  $\dD2$.
\end{thm}

\begin{proof}
  Let $f$ be a function which is infinite time eventually computable
  using the program $p$, and $(x,s)$ be given.  We can eventually
  decide whether $s\subset f(x)$ by simulating $p$ on input $x$ and
  checking at each stage whether $s$ is contained in the output.
\end{proof}

\begin{cor}
  \label{cor_meas}
  Every infinite time eventually decidable set is measurable.  Every
  infinite time eventually computable function is a measurable
  function.
\end{cor}

When we speak of measure, we are of course referring to the natural
coin-flipping probability measure on $2^\omega$, also called the
Lebesgue or Haar measure.  It is just the $\omega$-fold product of the
$(\frac12,\frac12)$ measure on $\set{0,1}$ The following result will
be of fundamental importance in later sections.

\begin{proof}
  By \cite[Exercise~14.4]{kanamori} every absolutely $\dD2$ set is
  measurable, and hence every infinite time eventually decidable set
  is measurable.  If $f$ is an infinite time eventually computable
  function, it follows from Theorem~\ref{thm_provd12} that for every
  open $U\subset2^\omega$, $f^{-1}(U)$ is absolutely $\dD2$.  Hence
  $f^{-1}(U)$ is measurable, and so $f$ is measurable.
\end{proof}

\section{Infinite time effective reductions}

We are now ready to begin our generalization of Borel reducibility to
the infinite time computability context. We shall introduce several
classical results of Borel equivalence relation theory in turn and
inquire how they are transferred or transformed to the infinite time
computability context. And we shall also introduce several of the
natural equivalence relations that we aim to fit into our new
hierarchy. In the following section, we will begin to treat the
infinite time analogue of the countable Borel equivalence
relations.

Recall that if $E,F$ are equivalence relations on $2^\omega$, then $f$
is a reduction from $E$ to $F$ if and only if it satisfies
\[x\mathrel{E}y\iff f(x)\mathrel{F}f(y)\;.
\]
We say that $E$ is Borel reducible to $F$, written $E\leq_BF$, if
there is a Borel reduction from $E$ to $F$. We propose to focus on the
following generalizations of the reduction concept to the context of
infinite time computability, corresponding to the two notions of
computability that we have discussed.
\begin{itemize}
\item The relation $E$ is infinite time \emph{computably reducible} to
  $F$, written $E\leq_cF$, if there is an infinite time computable
  reduction from $E$ to $F$.
\item The relation $E$ is infinite time \emph{eventually reducible} to
  $F$, written $E\leq_eF$, if there is an infinite time eventually
  computable reduction from $E$ to $F$.
\end{itemize}

To begin with some elementary considerations from the Borel theory,
let $\Delta(1)$, $\Delta(2),\ldots$, $\Delta(\omega)$ denote arbitrary
but fixed Borel equivalence relations with $1,2,\ldots,\omega$
classes, respectively.  Then
$\Delta(1)<_B\Delta(2)<_B\cdots<_B\Delta(\omega)$, and moreover these
are the simplest relations in the sense that for any $E$ with
infinitely many classes, $\Delta(\omega)\leq_BE$.  The next least
complex Borel equivalence relation is the \emph{equality relation} on
$2^\omega$, sometimes denoted $\Delta(2^\omega)$ or simply
$\mathord{=}$.

\begin{thm}[Silver dichotomy]
  If $E$ is a Borel (or even $\dP1$) equivalence relation then either
  $E$ has at most countably many classes or else
  $\,\mathord{=}\leq_BE$.
\end{thm}

Equivalence relations $E$ which are Borel reducible to $\mathord{=}$
are called \emph{smooth} or \emph{completely classifiable}, since the
corresponding reduction function shows how to concretely compute
complete invariants for $E$.  One step further up the hierarchy, one
finds the \emph{almost equality} relation $E_0$, which is defined by
$x\mathrel{E}_0y$ if and only if $x(n)=y(n)$ for almost all $n$.

We now present a proof that $E_0$ is not Borel reducible to
$\mathord{=}$.  This will be the first example of a proof that there
cannot be a Borel reduction from $E$ to $F$ which overshoots and shows
more.  In this case, it shows that there cannot be a measurable
reduction from $E$ to $F$, and hence there cannot be an infinite time
decidable or even provably $\dD2$ such reduction.  We shall discuss
this phenomenon further int he last section.

\begin{prop}
  \label{e0}
  There is no measurable reduction from $E_0$ to equality
  $\,\mathord{=}$, and hence $\,\mathord{=}<_cE_0$.
\end{prop}

\begin{proof}
  Suppose that $f$ is a measurable reduction from $E_0$ to
  $\mathord{=}$.  Then for every $U\subset2^\omega$, $f^{-1}(U)$ is
  closed under $E_0$ equivalence, \emph{i.e.} it is closed under
  finite modifications.  Such a set is called a ``tail set'', and a
  direct argument shows that such sets have measure $0$ or $1$.
  Letting $U$ run over the basic sets, we obtain that $f$ is constant
  on a set of measure $1$.  But $f$ is countable-to-one, and since the
  measure is nonatomic, this is a contradiction.
\end{proof}

We presently discuss a second dichotomy theorem (see \cite{hkl}).

\begin{thm}[Glimm-Effros dichotomy]
  \label{glimm-effros}
  If $E$ is any Borel equivalence relation then either $E$ is smooth
  or else $E_0\leq_BE$.
\end{thm}

Neither the Silver dichotomy nor the Glimm-Effros dichotomy can hold
in the case of infinite time decidable equivalence relations and
infinite time computable reductions for the simple reason that there
exist infinite time computable equivalence relations which necessarily
have $\aleph_1$ many classes.  But it is conceivable that this is the
only obstruction, and many questions about $E_0$ and infinite time
computable equivalence relations remain open.

\begin{question}
  Do any useful generalizations of the Silver dichotomy or
  Glimm-Effros dichotomy hold in the case of infinite time decidable
  equivalence relations and infinite time computable reductions?
\end{question}

One might ask if there is any difference whatsoever between the Borel
and infinite time computable theories. Of course not \emph{every}
infinite time computable reduction can be replaced by a Borel
reduction. To give a trivial counterexample, consider an infinite time
decidable equivalence relation $E$ with just two non-Borel classes:
clearly, in this case we have $E\leq_c\Delta(2)$ and
$E\not\leq_B\Delta(2)$.  Of course, such counterexamples must be
pervasive in the hierarchy, for instance by replacing one of those
equivalences classes with another entire equivalence relation.  In
Proposition \ref{wo_ck} we shall give a \emph{naturally occurring}
pair of equivalence relations such that $E\leq_cF$ and $E\not\leq_BF$.
However, our example will be of high descriptive complexity, and
so we are left with the following interesting problem.

\begin{question}
  Are there Borel equivalence relations $E,F$ such that $E\leq_cF$ but
  $E\not\leq_BF$?
\end{question}

For an example of such $E$ and $F$ of higher complexity, we consider
the following two equivalence relations.
\begin{itemize}
\item Let $x\iso_{WO} y$ if and only if $x$ and $y$, thought of as
  codes for binary relations on $\omega$, code isomorphic wellorders
  on $\omega$.
\item Let $x\mathrel{E}_{ck}y$ if and only if
  $\omega_{1,ck}^x=\omega_{1,ck}^y$, that is, if and only if $x$ and
  $y$ can write the same ordinals in $\omega$ steps.
\end{itemize}

\begin{prop}
  \label{wo_ck}
  The equivalence relations $E_{ck}$ and $\oiso_{WO}$ are infinite
  time computably bireducible.  On the other hand, they are Borel
  incomparable.
\end{prop}

\begin{proof}
  Results in \cite{machines} show that for any real $x$, the ordinal
  $\omega_{1,ck}^x$ is infinite time writable from parameter $x$, and
  this algorithm is uniform in $x$. So there is an infinite time
  computable function $f$ such that $f(x)$ is a real coding the
  ordinal $\omega_{1,ck}^x$. This function is therefore a reduction
  from $E_{ck}$ to $\oiso_{WO}$.  Next, we shall show there can be no
  Borel reduction $f$ from $E_{ck}$ to $\oiso_{WO}$.  Suppose that $f$
  is such a reduction.  It takes values in $WO$ and since $E_{ck}$ has
  $\omega_1$ many equivalence classes, the range of $f$ must code
  unboundedly many ordinals.  By the boundedness theorem (see
  \cite[Theorem 35.23]{kechris}), $\im(f)$ is not $\dS1$ and hence $f$
  is not Borel.

  Next we show that $\oiso_{WO}$ reduces to $E_{ck}$.  Let $y$ be a
  code for an ordinal $\alpha$.  We shall compute $x=f(y)$ depending
  only on $\alpha$ and such that $\omega_1^x$ is equal to the
  $\alpha^\mathrm{th}$ admissible ordinal $\delta$.  First, given a
  code $z$ for an ordinal $\beta$ we can always find its admissible
  successor (the least admissible above $\beta$).  To see this, note
  that it must be bounded by $\omega_{1,ck}^z$.  So for each
  $\beta<\gamma\leq\omega_{1,ck}^z$ we simply build $L_\gamma$ and
  check to see that it satisfies the \textsf{KP} axioms.  Now we can
  iterate this $\alpha$ times to find the $\alpha^\mathrm{th}$
  admissible ordinal $\delta$.  Next build $L_\delta$, and search
  inside it to find the $L$-least $x$ such that
  $\omega_{1,ck}^x=\delta$.  Clearly $x$ depends only on $\alpha$ and
  not the given code $y$.

  Finally, we argue that no such reduction $f$ can be Borel.  Indeed,
  notice that $E_{ck}$ is a $\dS1$ relation since $x\mathrel{E}_{ck}y$
  is equivalent to the following $\dS1$ assertion: Whenever $e$ is a
  finite time program such that $\varphi_e(x)$ codes a well order,
  there exists a finite time program $e'$ such that
  $\varphi_e(x)\cong\varphi_{e'}(y)$.  Hence, $\im(f)$ is a $\dS1$
  subset of $WO$.  Since $\im(f)$ is necessarily unbounded, this
  contradicts the boundedness theorem.
\end{proof}

Gao had noted \cite[Section~9.2]{gao} that there exists a $\dD2$
reduction from $\oiso_{WO}$ to $E_{ck}$, but that the study of $\dD2$
reducibility is problematic.  Thus Proposition~\ref{wo_ck} resolves
this by showing that the reduction from $\oiso_{WO}$ to $E_{ck}$ is
infinite time computable.  To see that the $\dD2$ reductions can pose
difficulties, we now show that in $L$ the $\dD2$ functions, and indeed
the semicomputable functions, collapse a large portion of the
hierarchy of equivalence relations.

\begin{thm}
  \label{thm_semi}
  If $V=L$ then whenever $E$ is an infinite time decidable equivalence
  relation, there exists an infinite time semicomputable function
  $f\from 2^\omega\rightarrow2^\omega$ such that $x\mathrel{E}y$ if
  and only if $f(x)=f(y)$.
\end{thm}

\begin{proof}
  Following the $L$-code argument of \cite[Theorem 38]{model}, given
  $x\in 2^\omega$ we shall encode its equivalence class by a pair of
  reals.  Let $\alpha<\omega_1$ be least such that $L_\alpha$ contains
  a member of the $E$-equivalence class of $x$, and $L_\alpha\models$
  ``some (fixed) large fragment of \ZFC\ and $\omega_1$ exists.''  The
  idea is that $L_\alpha$ is large enough that all computations on
  reals of $L_\alpha$ halt or repeat in fewer than $\alpha$ steps.

  Let $\beta>\alpha$ be least such that $\beta$ is countable in
  $L_{\beta+1}$ and let $w\in L_{\beta+1}$ be the $L$-least real
  coding $\beta$.  Finally let $z$ be the $L$-least real which is
  $E$-equivalent to $x$.  Then by our choice of $\alpha$, we have
  $z\in L_\alpha$.  Now, $f(x)\defeq w\oplus z$ is the code we
  seek.

  We clearly have $x\mathrel{E}y$ if and only if $f(x)=f(y)$, but we
  must verify that $f$ is infinite time semicomputable.  That is, a
  machine must recognize given $(x,w_0\oplus z_0)\in f$, whether
  $w_0=w$, $z_0=z$ as defined above.  The machine first checks to see
  that $w_0$ codes an ordinal, and using this ordinal as $\beta$ it
  constructs $L_\alpha$ and checks to see that $L_\alpha\models
  z_0=z$.  Lastly, note that $L_\alpha$ is correct about this since it
  has access to all computations on its reals.
\end{proof}

On the other hand, we have seen that infinite time effective sets and
functions derive many of their properties from the fact that they are
\emph{absolutely} $\dD2$.  It is therefore natural to study absolutely
$\dD2$ reducibility, as we shall do in the last section.  One might
therefore ask whether there is any sense in which the absolutely
$\dD2$ sets and functions are effective.  The following result sheds
doubt on this by showing first that this question cannot be separated
from that of whether there is a sense in which \emph{all} $\dD2$ sets
are effective.

\begin{prop}
  There is a forcing extension of the universe in which every $\dD2$
  set is absolutely $\dD2$, and indeed, in which every equivalent pair
  of $\Sigma^1_2$ and $\Pi^1_2$ definitions remains equivalent after
  any further forcing.
\end{prop}

\begin{proof}
  Let us first show that the desired situation holds under the
  Maximality Principle, which is the scheme asserting that any
  forcibly necessary set-theoretic statement is already true (see
  \cite{Hamkins2003:MaximalityPrinciple}, also
  \cite{StaviVaananen:ReflectionPrinciples}). A statement is forceably
  necessary, if it is forceable in such a way that it remains true in
  all further forcing extensions. If $V$ satisfies the Maximality
  Principle and $\varphi$ and $\psi$ are $\Sigma^1_2$ and $\Pi^1_2$
  assertions, respectively, which could become inequivalent in a
  forcing extension $V[G]$, then there is a real $z$ in $V[G]$ such
  that $\varphi(z)$ differs from $\psi(z)$ in $V[G]$. Since these
  statements are each absolute to all further extensions of $V[G]$,
  this means that the inequivalence of $\varphi$ and $\psi$ is
  forceably necessary over $V$ and therefore true there by the
  Maximality Principle. Thus, under the Maximality Principle, any two
  $\Sigma^1_2$ and $\Pi^1_2$ assertions that are equivalent in $V$
  remain equivalent in all forcing extensions. In particular, every
  $\dD2$ set in $V$ is absolutely $\dD2$.

  This argument makes use of only a small fragment of the Maximality
  Principle. And although it is proved in
  \cite{Hamkins2003:MaximalityPrinciple} that if \ZFC\ is consistent,
  then there is a model of \ZFC\ plus the Maximality Principle, it is
  also observed there that some models of \ZFC\ have no forcing
  extensions with the Maximality Principle. Nevertheless, the main
  argument of \cite{Hamkins2003:MaximalityPrinciple} does show that
  every model of \ZFC\ has a forcing extension with the Maximality
  Principle restricted to assertions of a given fixed set-theoretic
  complexity. Since we only used low projective complexity in the
  previous paragraph, there is a forcing extension of the universe in
  which every $\dD2$ set is absolutely $\dD2$ as described. (The
  forcing is simply an iteration, where one continues forcing until
  all possible inequivalences have been exhibited.)
%
\end{proof}

On the other hand, there are models with $\dD2$ functions which are
not absolutely $\dD2$.  For instance in $L$ there is a $\dD2$
well-ordering of the reals, though no model of \ZFC\ has an absolutely
$\dD2$ well-ordering of the reals.

We close this section by introducing a number of equivalence relations
which are of natural interest.

\begin{itemize}
\item Let $x\mathrel{E}_{\textrm{set}}y$ if and only if $x$ and $y$,
  thought of as countable sequences of reals, have the same range.
\item Let $x\iso_{HC}y$ if and only if $x$ and $y$, thought of as
  codes for hereditarily countable sets, are isomorphic. (Here, $x$ is
  said to code a hereditarily countable set $z$ iff, thinking of $x$
  as a binary relation on $\omega$, we have
  $(\omega,x)\iso(tc(\set{z}),\mathord{\in})$.)
\item Let $x\iso y$ if and only if $x$ and $y$, thought of as codes
  for countable structures in a countable language, are isomorphic.
\item Let $x\mathrel{E}_\lambda y$ if and only if
  $\lambda^x=\lambda^y$, that is, if and only if $x$ and $y$ can write
  the same set of ordinals.  Similarly, define $x\mathrel{E}_\zeta y$
  if and only if $x$ and $y$ can eventually write the same ordinals,
  and $x\mathrel{E}_\Sigma y$ if and only if $x$ and $y$ can
  accidentally write the same ordinals.
\item Let $x\equiv_Ty$ if and only if $x,y$ lie in the same classical
  Turing degree.
\item Let $x\equiv_\textrm{arith}$ if and only if $x,y$ lie in the
  same arithmetic degree.
\item Let $x\equiv_\textrm{hyp}$ if and only if $x,y$ lie in the same
  hyperarithmetic degree.
\item Let $x\equiv_\infty y$ if and only if $x$ and $y$ are infinite
  time computable from one another (that is, lie in the same infinite
  time \emph{degree}).
\item Let $x\equiv_{e\infty} y$ if and only if $x$ and $y$ are
  infinite time eventually computable from one another (that is, lie
  in the same infinite time \emph{eventual degree}).
\item Let $x\mathrel{J}y$ if and only if $x,y$ have equivalent
  infinite time jumps, \emph{i.e.}, $x^\triangledown\equiv_\infty
  y^\triangledown$.
\end{itemize}

Notice that $x\iso_{WO}y$ only makes sense for those $x,y\in2^\omega$
which code a well order. Thus, this is a relation not on all of Cantor
space, but only on the (infinite time decidable) subset consisting of
codes for well orders.  This issue of an equivalence relation that is
merely partial never arises in Borel equivalence relations since the
domain of any Borel equivalence relation is a standard Borel space in
its own right, and so every Borel equivalence relation can be assumed
to be total.  We shall allow for the study of infinite time computable
relations $E\subset 2^\omega\times2^\omega$ whose domain is defined on
an infinite time computable subset of $2^\omega$.

Some reductions between these equivalence relations are already
apparent.  For instance, given a countable sequence of reals $\langle
a_n\rangle$, it is not difficult to construct an HC-code for the set
$\{a_n\}$, and hence $E_{\textrm{set}}$ is computably reducible to
$\oiso_{HC}$.  Next, $\oiso_{WO}$ is computably reducible to
$\oiso_{HC}$ since it is just the restriction of this relation to the
set $WO$ of codes for well orders.  Thirdly, it is easy to see that
the function $x\mapsto x^\triangledown$ is an eventual reduction from
$J$ to $\oequiv_{\infty}$.  Many more details of the
interrelationships (with respect to infinite time computable and
infinite time eventually computable reducibility) shall be examined as
the exposition unfolds.

\section{Enumerable equivalence relations}

The classical Borel equivalence relation theory has placed a major
focus on the countable Borel equivalence relations, and the
investigation of this natural sub-hierarchy of the hierarchy of all
equivalence relations has led to some of the most fruitful work (see
for instance \cite{jkl}).  Not only do these relations include many of
the most natural examples, but some of most powerful methods in the
theory appear to work most effectively with countable relations. The
situation is rather reminiscent of the focus in computability theory
on the c.e.\ Turing degrees as a sub-hierarchy of the hierarchy of all
Turing degrees.

An equivalence relation $E$ is \emph{countable} if every
$E$-equivalence class is countable.  A key characterization of the
countable Borel equivalence relations is that they are exactly those
relations $E$ with a \emph{Borel enumeration}, a Borel function $f$
such that $f(x)=\langle x_0,x_1,\ldots\rangle$ effectively enumerates
the elements of $[x]_E$. (This is a consequence of the Lusin-Novikov
theorem, \cite[Theorem~18.10]{kechris}.)  The natural extension of the
class of countable Borel equivalence relations to the infinite time
computable context simply generalizes this enumerability concept.

\begin{defn}\ 
  \begin{itemize}
  \item A countable equivalence relation $E$ is infinite time
    \emph{enumerable} if it admits an infinite time computable
    enumeration function, that is, a function $f$ for which
    $f(x)=\langle x_0,x_1,\ldots\rangle$ enumerates $[x]_E$ for all
    $x\in2^\omega$.
  \item Similarly, $E$ is infinite time \emph{eventually enumerable}
    if it admits an infinite time eventually computable enumeration
    function.
  \end{itemize}
\end{defn}

Recall that if $\Gamma$ is any group of bijections of $2^\omega$, we
can define the corresponding \emph{orbit equivalence relation}
$E_\Gamma$ by
\[x\mathrel{E}_\Gamma y\iff\Gamma x=\Gamma y.
\]
By a theorem of Feldman and Moore \cite{feldmanmoore}, $E$ is a
countable Borel equivalence relation if and only if there exists a
countable group $\Gamma$ of Borel bijections of $2^\omega$ such that
$E=E_\Gamma$.  Our first observation is that the infinite time
enumerable relations enjoy an analogous property.

\begin{thm}
  \label{thm_enumerable_char}
  An equivalence relation $E$ is infinite time enumerable if and only
  if there exists a countable group $\Gamma$ of infinite time
  computable bijections of $2^\omega$ such that $E$ is precisely the
  induced orbit equivalence relation $E_\Gamma$.  The analogous result
  holds for the infinite time eventually enumerable equivalence
  relations.
\end{thm}

\begin{proof}
  Suppose first that there exists such a group $\Gamma$.  Write
  $\Gamma=\langle\gamma_n\rangle$, and let $r$ be a real code for a
  sequence $\langle r_n\rangle$ such that each $r_n$ codes a program
  which computes the function $\gamma_n$.  We claim that $E_\Gamma$ is
  infinite time enumerable in the real $r$.  Indeed, on input $x$, a
  program can simply use $r$ to simulate each $\gamma_n$ on input $x$,
  and collect the values $\gamma_n(x)$ into a sequence.

  Conversely, suppose that $E$ is infinite time enumerable.  By the
  proof of the classical Feldman-Moore theorem, it suffices to
  establish the conclusion of the Lusin-Novikov theorem, namely:
  \begin{itemize}
  \item $E$ can be expressed as a countable union of graphs of
    infinite time computable partial functions.
  \end{itemize}
  For this, let $f$ be an infinite time computable function which
  witnesses that $E$ is infinite time enumerable, \emph{i.e.}, for
  every $x\in2^\omega$, $f(x)$ is a code for the $E$-class of $x$.
  Letting $f_n(x)$ denote the $n^\mathrm{th}$ element of $f(x)$, we
  have that $E=\cup f_n$.  This completes the proof.
\end{proof}

\begin{prop}
  The class of infinite time enumerable equivalence relations lies
  properly between the countable Borel equivalence relations and the
  countable infinite time decidable equivalence relations.
\end{prop}

\begin{proof}
  That every countable Borel equivalence relation is infinite time
  enumerable follows from the previous theorem, and it is immediate
  from the definition that every infinite time enumerable equivalence
  relation is countable and infinite time decidable.

  We now give an example of a countable infinite time decidable
  equivalence relation which is not infinite time enumerable.  For
  each $x\in2^\omega$, we let $c^x$ denote the lost melody real
  relative to $x$.  Recall that $c^x$ is a real such that $\{c^x\}$ is
  (lightface) infinite time decidable in $x$ and yet $c^x$ is not
  infinite time writable in $x$.  It follows that the function
  $f(x):=x\oplus c^x$ is infinite time semicomputable but not infinite
  time computable, even from a real parameter.  Now, we let
  $x\mathrel{E}y$ if and only if there exists $n$ such that $x=f^n(y)$
  or $y=f^n(x)$.  Since $f$ is injective, $E$ is an equivalence
  relation. Moreover, it is easy to see that $E$ is countable and
  infinite time decidable.  However, $E$ cannot be infinite time
  enumerable in the parameter $z$, for then $c^z$ would be infinite
  time writable in $z$, a contradiction. Indeed, $E$ cannot even be
  accidentally enumerable.

  For an example of an infinite time enumerable equivalence relation
  which is not Borel, we shall use hyperarithmetic equivalence
  $\oequiv_\textrm{hyp}$.  Recall that $x\equiv_\textrm{hyp}y$ if and
  only if $x\in\Delta_1^1(y)$ and $y\in\Delta_1^1(x)$.  It follows
  from the proof of Theorem~\ref{thm_Borel_char} that
  $x\equiv_\textrm{hyp}y$ if and only if $x$ is infinite time
  computable from $y$ in fewer than $\omega_{1,ck}^y$ steps and $y$ is
  infinite time computable from $x$ in fewer than $\omega_{1,ck}^x$
  steps.  (Recall that $\omega_1^y$ denotes the supremum of the
  ordinals computable in the ordinary sense from $y$.)

  Since $\omega_{1,ck}^x$ is infinite time computable from $x$, the
  equivalence relation $\oequiv_\textrm{hyp}$ is clearly infinite time
  enumerable.  But suppose, towards a contradiction, that
  $\oequiv_\textrm{hyp}$ is Borel.  Then since $\oequiv_\textrm{hyp}$
  is also countable, there exists a Borel function $f$ such that for
  all $x$, $f(x)$ codes $[x]_{\oequiv_\textrm{hyp}}$.  By Theorem
  \ref{thm_Borel_char}, there exists a program $e$ in a parameter $z$
  and an ordinal $\alpha$ such that on any input $x$, $e$ computes
  $f(x)$ in fewer than $\alpha$ steps.  Replacing $z$ with a more
  complicated real if necessary, we may suppose that
  $\alpha\leq\omega_{1,ck}^z$.  Now, using $e$ it is easy to write a
  program which first enumerates $[z]_{\oequiv_\textrm{hyp}}$, then
  diagonalizes against this set to write a real $r=z\oplus d$ such
  that $r\not\in[z]_{\oequiv_\textrm{hyp}}$.  Since $z$ is quickly
  writeable from $r$, we must have that $r$ isn't writable from $z$ in
  fewer than $\omega_{1,ck}^z$ steps.  This is a contradiction,
  because we have just described a program which does so.
\end{proof}

The infinite time eventually enumerable equivalence relations are
easily seen to be infinite time semidecidable, but as the next
proposition shows, not necessarily infinite time enumerable or even
infinite time decidable.

\begin{prop}
  $\oequiv_\infty$ is infinite time eventually enumerable but not
  infinite time decidable.
\end{prop}

\begin{proof}
  That $\oequiv_\infty$ is infinite time semidecidable is shown in
  \cite[Theorem 5.7]{machines}; the argument is very simple.  On input
  $x,y$, just simulate all programs on input $x$ and see if any of
  them writes $y$, and vice versa.

  Now, suppose towards a contradiction that $\oequiv_\infty$ is
  infinite time decidable in the real parameter $z$.  We shall use
  this to decide the halting problem in $z$, \emph{i.e.} the real
  $z^\triangledown=\{e:\varphi_e^z(0)\;\mathrm{halts}\}$.  Consider
  the program which attempts to compute this set.  It runs all
  programs simultaneously on input $0$, and each time one halts, its
  output is added to an accumulating set $x$.  Additionally, it checks
  at each stage whether $x\oplus z\equiv_\infty z$ and halts if this
  does not hold.

  Note that this program halts, since after some stage all programs
  which halt have halted.  Moreover, at this moment the approximation
  is correct and so certainly $x\oplus z\not\equiv_\infty z$.  It may
  halt earlier than this, but it must halt with some real $x$ such
  that $x\oplus z\not\equiv_\infty z$.  Hence from $z$ it has computed
  a real strictly more complex than $z$, a contradiction.
\end{proof}

In particular, we have the following consequence.

\begin{cor}
  The relation $\oequiv_\infty$ doesn't computably reduce to any
  infinite time decidable equivalence relation.
\end{cor}

\begin{prop}
  The relation $\oequiv_{e\infty}$ is infinite time eventually
  decidable but not infinite time eventually enumerable.
\end{prop}

That the infinite time eventual degree relation $\oequiv_{e\infty}$ is
infinite time eventually decidable is due to Philip Welch.  This
result was very surprising to experts in the area, since the
corresponding infinite time Turing degree relation $\equiv_{\infty}$
was known not be be infinite time decidable in \cite{machines} by the
elementary argument above.

\begin{proof}
  First, to see that the set of infinite time eventually writable
  reals is not infinite time eventually enumerable, suppose that $e$
  is a program with oracle $z$ which on input $x$, eventually writes a
  code for the $\oequiv_{e\infty}$-class of $x$.  Then consider the
  diagonalization program $q$ which simulates $e$ on $z$, and at each
  stage of simulation writes a real which is not in the set coded on
  the output tape of $e$. Then $q$ eventually writes a real which is
  not in the $\oequiv_{e\infty}$-class of $z$, a contradiction.

  Now, to show that $\oequiv_{e\infty}$ is infinite time eventually
  decidable, we shall actually just show that the set of infinite time
  eventually writable reals is infinite time eventually decidable.
  The proposition follows, since this argument relativizes to show
  that given $z$, the set of reals infinite time eventually writable
  in $z$ is infinite time eventually decidable.  Following the
  infinite time Turing machine literature, we denote by $\lambda$ the
  supremum of the writable ordinals (which is the same as the supremum
  of the clockable ordinals), by $\zeta$ the supremum of the infinite
  time eventually writable ordinals and by $\Sigma$ the supremum of
  the infinite time accidentally writable ordinals. Results in
  \cite{machines} establish that $\lambda<\zeta<\Sigma$, and Philip
  Welch (\cite{length}, \cite{Welch2000:Eventually}, see also
  \cite[Theorem~1.1]{post}) has proved moreover that
  $L_\lambda\prec_{\Sigma_1} L_\zeta \prec_{\Sigma_2} L_\Sigma$, and
  furthermore these ordinals are characterized as least having this
  property.  This key result is now known as the
  $\lambda$-$\zeta$-$\Sigma$ theorem. Welch proved that every infinite
  time Turing machine computation either halts in time before
  $\lambda$ or repeats its stage $\zeta$ configuration at
  $\Sigma$. Any computation that eventually stabilizes, reaches its
  stabilizing configuration before $\zeta$, and the universal
  computation simulating all programs on trivial input repeats the
  stage $\zeta$ configuration at stage $\Sigma$ for the first
  time. Because of this, it is infinite time decidable whether a given
  real codes the ordinal $\Sigma$, since the machine need merely check
  that it does indeed code a well order, and that the universal
  computation, when simulated for that many steps, exhibits this
  repeating phenomenon exactly at that stage. These facts relativize
  easily to a real parameter.

  Now, on input $x$, we can eventually decide whether it is infinite
  time eventually writable by the following algorithm.  First, write a
  preliminary default ``No'' on the output tape. Next, simulate the
  universal computation, and search to see if $x$ is ever shown to be
  accidentally writable. If so, change the answer on the output tape
  provisionally to ``Yes,'' and then run the universal program with
  parameter $x$ to see if there is an $x$-writable real coding the
  ordinal $\Sigma$. By the remarks in the previous paragraph, any
  instance of this is infinite time decidable. If $\Sigma$ is ever
  found to be $x$-writable, then change the answer finally back to
  ``No'' and halt; otherwise, keep searching.

  Let's argue that this algorithm works. If $x$ is eventually
  writable, then it will appear on the tape before stage $\zeta$, and
  so we shall pass the first hurdle, where the answer was changed
  provisionally to ``Yes.'' But since $x$ is eventually writable, it
  follows that $\Sigma$ cannot be $x$-writable, since if it were, then
  $\Sigma$ would be accidentally writable, contradicting the fact that
  is larger than all accidentally writable ordinals. Thus, in this
  case we shall never pass the second hurdle, and so our final answer
  will stabilize on ``Yes,'' as desired. If $x$ is not eventually
  writable and also not accidentally writable, then the algorithm will
  never pass the first hurdle, and so the algorithm will stabilize on
  the first ``No,'' as desired. Finally, if $x$ is accidentally
  writable, but not eventually writable, then $x$ appears accidentally
  on the universal computation, but not before time $\zeta$ (since
  otherwise it would be eventually writable). So it appears at some
  point between $\zeta$ and $\Sigma$. In this case, the ordinal
  $\zeta$ is below the supremum of the $x$-clockable ordinals, and
  since the supremum of the $x$-clockable and $x$-writable ordinals is
  the same, it follows that $\zeta$ is $x$-writable. From this, it
  follows that there are ordinals above $\Sigma$ that are
  $x$-clockable, since with oracle $x$ we can run the universal
  computation, look exactly at the stage $\zeta$ configuration, and
  then wait until stage $\Sigma$, when this configuration first
  repeats.  Thus, $\Sigma$ is also $x$-writable, and so the algorithm
  will pass the final hurdle, changing the answer to ``No,'' and
  halting, as desired.
\end{proof}

\begin{figure}
\begin{center}
\includegraphics[width=2in]{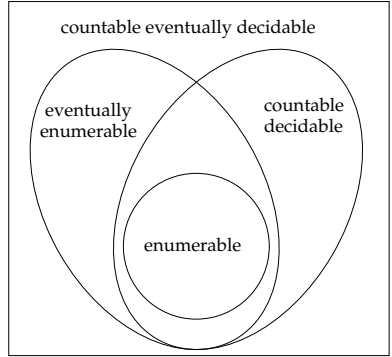}
\end{center}
\caption{The relationships between various classes of countable
  equivalence relations.\label{fig_enumerable}}
\end{figure}

The relationships between various classes of countable equivalence
relations are shown in Figure~\ref{fig_enumerable}.  Each of the
inclusions is proper; we have only omitted the fact that there exists
an infinite time eventually enumerable, infinite time decidable
equivalence relation which is not infinite time enumerable.

We now turn towards an analysis of the structure of the infinite time
enumerable equivalence relations.  We begin by describing the most
basic structure theory of the countable Borel equivalence
relations. First, we have already seen that it is a consequence of
Silver's theorem that the equality relation $\mathord{=}$ is the
minimum countable Borel equivalence relation.  The relations $E$ which
are Borel reducible to $\mathord{=}$ are called \emph{smooth}.  By the
Glimm-Effros dichotomy (Theorem \ref{glimm-effros}), $E_0$ is the
next-least countable Borel equivalence relation, in the strong sense
that $E_0$ is Borel reducible to \emph{any} nonsmooth Borel
equivalence relation.  Lastly, and somewhat surprisingly, there exists
a \emph{universal} countable Borel equivalence relation, denoted
$E_\infty$.  It is realized as the orbit equivalence relation induced
by the left-translation action of the free group $F_2$ on its power
set.

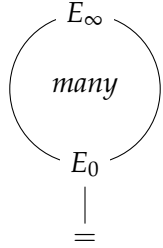
\begin{figure}
\begin{center}
\begin{tikzpicture}
\draw (0,2) circle (1) ;
\node at (0,0) (equality) {$=$} ;
\node[fill=white] at (0,1) (e0) {$E_0$}
  edge (equality) ;
\node at (0,2) {\emph{many}} ;
\node[fill=white] at (0,3) (einfty) {$E_\infty$} ;
\end{tikzpicture}
\end{center}
\caption{The countable Borel equivalence relations\label{fig_cbers}.}
\end{figure}

There were initially very few countable Borel equivalence relations
known to lie in the interval $(E_0,E_\infty)$.  It is a fundamental
result of Adams and Kechris \cite[Theorem 1]{adamskechris} that there
exists a sequence $\{AK_\alpha\}$ of continuum many pairwise
incomparable countable Borel equivalence relations.  In summary, we
have that the countable Borel equivalence relations are organized as
in Figure \ref{fig_cbers}.

We would like to develop an analogous picture for the infinite time
enumerable relations.  We first consider the question of whether the
Silver dichotomy holds for the infinite time enumerable relations,
that is, whether $\mathord{=}$ is the least complex such relation.

\begin{thm}\ 
  \label{thm_silver}
  \begin{itemize}
  \item There is a perfect set of infinite time eventual degrees.
  \item If $E$ is infinite time eventually enumerable then $\mathord{=}\leq_BE$.
  \end{itemize}
\end{thm}

The proof hinges on the following result.  It is nearly implicit in
\cite{minimality}, and Welch has subsequently completed the proof
based on that work.  We shall shortly provide a different argument
which is due to Hamkins.

\begin{thm}[Welch]
  \label{thm_generic}
  If $c$ is an $L_{\Sigma}$-generic Cohen real, then
  $\lambda^c=\lambda$, $\zeta^c=\zeta$ and $\Sigma^c=\Sigma$.
\end{thm}

We remark that our proof of this result will easily relativize.  That
is, for any real $z$, if $c$ is an $L_{\Sigma^z}[z]$-generic Cohen
real, then we shall have $\lambda^{z+c}=\lambda^z$,
$\zeta^{z+c}=\zeta^z$ and $\Sigma^{z+c}=\Sigma^z$.  Admitting this
result, let us show how to complete the proof of
Theorem~\ref{thm_silver}.

\begin{proof}[Proof of Theorem~\ref{thm_silver}]
  We begin by arguing that there is a perfect set of eventual degrees.
  Since $L_\Sigma$ is countable, there exists a perfect set $\mathcal
  G$ of reals which are mutually generic over $L_\Sigma$.  It suffices
  to show that for $g,g'\in\mathcal G$, if $g\neq g'$ then
  $g\not\equiv_{e\infty} g'$.  Indeed, since $g,g'$ are mutually
  generic we have $g'\not\in L_\lambda[g]$.  Furthermore, since
  $\zeta^{g}=\zeta$, it follows that $g'\not\in L_{\zeta^g}[g]$, in
  other words $g'$ is not infinite time eventually writable from $g$.

  Next, let $x\equiv_{e\infty}^zy$ if and only if $x$ and $y$ are
  infinite time eventually writable from one another using the
  parameter $z$.  Then by our earlier remarks,
  Theorem~\ref{thm_generic} relativizes and we in fact have
  $\mathord{=}\leq_B\oequiv_{e\infty}^z$ for any $z\in2^\omega$.
  Clearly, if $E$ is infinite time eventually enumerable in the
  parameter $z$ then $E\subset\oequiv_{e\infty}^z$, and it follows
  that $\mathord{=}\leq_eE$ as well.  The result for infinite time
  computable reducibility can be argued similarly.
\end{proof}

It is worth remarking that this argument also gives a reduction (the
same function) from $\mathord{=}$ to $J$.  We now return to the proof
of Theorem~\ref{thm_generic}.  Recall the result of Welch we mentioned
earlier, that for any $z$, every computation in $z$ either halts
before $\lambda^z$ or repeats the stage $\zeta^z$ configuration by
stage $\Sigma^z$, and moreover that the universal computation in $z$
repeats for the first time with this pair of ordinals.

\begin{proof}[Proof of Theorem~\ref{thm_generic}]
  Our strategy will be to show that for any infinite time Turing
  machine program $e$, $\varphi_e(c)$ repeats from stage $\zeta$ to
  stage $\Sigma$. Applying this to the case when $e$ is the universal
  program, this implies that $\zeta^c=\zeta$ and
  $\Sigma^c=\Sigma$. After this, we shall argue separately that
  $\lambda^c=\lambda$.

  The main idea is that instead of carrying out the computation
  $\varphi_e(c)$, which only exists in a world having $c$, we shall
  instead carry out a Boolean-valued computation using only the
  canonical name $\dot c$ for the Cohen generic, which is coded by a
  real in the ground model. The inspiration here is that if $c$ is
  fully $V$-generic, then every fact or aspect about the computation
  $\varphi_e(c)$, whether a given cell shows a $0$ or $1$ at a
  particular ordinal stage or whether the head is on a particular cell
  at a particular stage, is forced by some finite piece of the generic
  real $c$. This is the magic of forcing. We shall simply design a
  computation that keeps careful track of this information.

  Let us now describe the Boolean computation or simulation of
  $\varphi_e(\dot c)$. We embark on a computation that simulates
  $\varphi_e(\dot c)$ by computing exactly what information about this
  computation is forced by which conditions. At each simulated stage
  of computation, the algorithm keeps track of the values of the
  cells, the head position and the machine state, not with certainty,
  but with its corresponding Boolean value. That is, for each cell in
  the simulated computation, we reserve space in our actual
  computation to keep track of the conditions $p\in\PP$ that force the
  current value of this cell to be $0$ or to be $1$. Similarly, we
  also keep track of the conditions that force that the head is
  currently at this cell, and for each state in the program $e$ we
  keep track of the conditions that force that the simulated machine
  is currently in that state.

  Initially, our simulation data should specify that all conditions
  force that the head is on the left-most cell and in the
  \textsf{start} state, and that all the cells on the work and output
  tapes are $0$. For the input tape, which we intend to hold the
  generic real $c$, for each cell $j$ we say that a condition $p$
  forces that the $j\th$ cell has value $p(j)$, if this is
  defined. Every condition forces that the cells on the scratch and
  output tapes are all initially $0$.

  At successor stages of simulation, we can easily update this data so
  as to carry out the simulated computational step.  For example, if
  $p$ forces that the head is at a certain position, reading a certain
  value and in a certain state, then we can adjust our data for the
  next step so that $p$ forces the appropriate values and head
  position after one step of the program $e$. There is a subtle
  tidying-up issue, in that it could happen at a successor step that
  after this update, although previously a condition $q$ did not
  force, say, a certain head position, nevertheless now $q\concat0$
  and $q\concat 1$ both force the same head positions (perhaps having
  arrived from different directions). In this case, we would want to
  say that $q$ also forces this head position. More generally, if the
  collection of conditions forcing a particular feature (cell value,
  head position or state) is dense below a given condition $q$, then
  in our update procedure we tidy-up our data to show that $q$ also
  forces this feature.

  Let us now explain how to update the data at limit stages of
  computation. Of course, at any simulated limit stage, we want every
  condition to force that the head is now on the left-most cell in the
  \textsf{limit} state. It is somewhat more subtle, however, to update
  the cells on the tape correctly. The problem is that the data we now
  have available is the $\limsup$ of the previous data for the cell
  values, which is is not the same as the data for the $\limsup$ of
  the cell values. Nevertheless, we will be able to recover the data
  we need. Note that $p$ forces a particular cell value is $0$ at
  limit time $\alpha$ if and only if there are densely many $q$ below
  $p$ such that for some $\beta<\alpha$, the condition $q$ forces that
  this cell is $0$ from $\beta$ up to $\alpha$. But we do have this
  information available in the $\limsup$ of the previous data, since
  if $q$ forced the cell was $0$ from $\beta$ up to $\alpha$, then the
  limit of this data will continue to show that. Thus, we can
  correctly compute the correct forcing relation for the cell values
  on the tape at limit stages of the simulated computation. A simple
  inductive argument on the length of the computation now establishes
  that we have correctly calculated the forcing relation for the head
  position, machine states and cell values at every stage of simulated
  computation.

  If $c$ is actually generic, then this Boolean-valued computation
  collapses to the actual computation of $\varphi_e(c)$ as follows. At
  every stage $\alpha$, there are dense sets of conditions $p$ forcing
  exactly where the head is, and what the state is, and what appears
  in each particular cell. By genericity, the generic filter will meet
  each of these dense sets, and so as far as conditions in $c$ are
  concerned, the ghostly Boolean-valued computation follows along with
  the actual computation $\varphi_e(c)$. Indeed, we claim that if $c$
  is merely $L_\Sigma$-generic (meaning that $c$ meets all predense
  sets for $\PP$ coded as elements of $L_\Sigma$), then the
  computation of $\varphi_e(c)$ at stage $\alpha$ is exactly what is
  forced by some condition $p\in c$ in the Boolean-valued
  computation. This is certainly correct at the initial stage, and by
  induction it is preserved through successor stages and limits,
  because all the relevant dense sets are in $L_\Sigma$, and so $c$
  meets them as required in order to collapse the Boolean-valued
  computation.

  The key observation now is that since the Boolean-valued computation
  is repeating from $\zeta$ to $\Sigma$, it follows that the true
  computation $\varphi_e(c)$ must also be repeating from $\zeta$ to
  $\Sigma$. And this is precisely what we had set out to prove. So we
  have established that $\zeta^c=\zeta$ and $\Sigma^c=\Sigma$.

  We finally argue that $\lambda^c=\lambda$. Suppose that some $e$ had
  the property that $\varphi_e(c)$ halts at some ordinal stage
  $\alpha>\lambda$. Then some condition $p$ forces that $\varphi_e(c)$
  halts at stage $\alpha$. We may now run Boolean-valued computation
  and wait until $p$ forces that \textsf{halt} is achieved. Since the
  simulated computation takes at least as long as the actual
  computation, this would allow us to halt beyond $\lambda$, a
  contradiction, since there are no clockable ordinals above
  $\lambda$.
\end{proof}

A slew of questions follows.  For instance, we have just seen in
Theorem~\ref{thm_silver} that there is a perfect set of eventual
degrees, and hence of infinite time computable degrees.  It is natural
to ask just how complex the infinite time Turing degree relations
$\oequiv_\infty$ and $\oequiv_{e\infty}$ actually are.

\begin{question}
  \label{question_e0}
  Does $E_0$ reduce (in any reasonable sense) to either
  $\oequiv_\infty$ or $\oequiv_{e\infty}$?
\end{question}

We have also just seen that a Silver dichotomy holds for infinite time
enumerable relations.  This leaves open the following related
question.

\begin{question}
  \label{question_glimmeffros} Does a Glimm-Effros dichotomy hold for
the infinite time enumerable equivalence relations?  In other words,
for any infinite time enumerable equivalence relation $E$ do we have
either $E\leq_c\Delta(2^\omega)$ or $E_0\leq_cE$?  (And similarly for
infinite time eventually enumerable relations with respect to eventual
reducibility.)
\end{question}

We next address the question of whether there is a universal infinite
time enumerable equivalence relation.

\newpage

\begin{thm}\ 
  \label{thm_universal_enumerable}
  \begin{itemize}
  \item If $E$ is infinite time enumerable, then $E\leq_cE_\infty$.
  \item If $E$ is infinite time eventually enumerable then $E\leq_eE_\infty$.
  \end{itemize}
\end{thm}

\begin{proof}[Sketch of proof]
  This is analogous to the proof that any countable Borel equivalence
  relation is Borel reducible to $E_\infty$.  In that argument, the
  key point is that any countable Borel equivalence relation can be
  expressed as the orbit equivalence relation induced by the Borel
  action of a countable group.  For our result, the key point is
  Theorem \ref{thm_enumerable_char}.
\end{proof}

In particular, $\oequiv_\infty$ is eventually reducible to
$E_\infty$. It is now natural to extend Question~\ref{question_e0} to
the following stronger statement.

\begin{question}
  \label{question_univ}
  Does $E_\infty$ reduce (in any sense) to $\oequiv_\infty$?  In other
  words, is $\oequiv_\infty$ universal infinite time eventually
  enumerable?
\end{question}

Note that Slaman has shown that $\oequiv_\mathrm{arith}$ is universal
countable Borel.  On the other hand, it us unknown whether $\oequiv_T$
is universal countable Borel.  If it is, then the Martin Conjecture
must fail.  For a discussion of this question see \cite{thomas-mc}.

Finally, we consider the question of whether there are incomparable
infinite time enumerable relations.  Indeed, it is not difficult to
see from the proof of the Adams-Kechris theorem that there cannot even
be a measurable reduction between any two $AK_\alpha$.  Hence, we
obtain the following result for free.

\begin{thm}[Adams-Kechris]
  There is a sequence $\{AK_\alpha\}$ of continuum many infinite time
  enumerable equivalence relations which are pairwise infinite time
  (eventually) computably incomparable.
\end{thm}

We conclude this section with a question regarding the following chain
of refinements of $\oequiv_\infty$.

\begin{defn}
  For $\alpha<\omega_1$, let $x\equiv_\alpha y$ if and only if $x$ and
  $y$ are infinite time computable from one one another (without
  parameters) by computations which halt in fewer than $\alpha$ steps.
  This is an equivalence relation whenever $\alpha$ is additively
  closed.
\end{defn}

\begin{prop}
  Each equivalence relation $\oequiv_\alpha$ is countable and Borel.
\end{prop}

\begin{proof}
  The main point of interest is that $\oequiv_\alpha$ are Borel.
  First note that $\oequiv_\alpha$ is infinite time computable from an
  oracle for a real coding $\alpha$, and it is easily seen that it is
  infinite time computable uniformly in at most $\alpha+\alpha$ steps.
  But by Theorem \ref{thm_Borel_char}, any uniformly infinite time
  decidable set is Borel.
\end{proof}

We have the equalities $\oequiv_0$ is $\Delta(2^\omega)$ and
$\oequiv_\omega$ is $E_0$.  The union of the $\oequiv_\alpha$ is again
countable, it is precisely $\oequiv_\infty$.  Nothing is known about
the rest of the $\oequiv_\alpha$ for $\alpha<\omega_1$.

\begin{question}
  What is the structure of the $\oequiv_\alpha$ under infinite time
  computable comparability?  Are they linearly ordered with respect to
  infinite time computable reducibility?
\end{question}

\section{Some tools for showing non-reducibility}

In this section we will establish several non-reducibility results,
that is, results which state that some equivalence relation $E$ is not
reducible to another equivalence relation $F$.  As we have mentioned,
many such non-reducibility results from the theory of Borel
equivalence relations come from arguments which shows that there
cannot be an absolutely $\dD2$ reduction from $E$ to $F$.  In this
section, we give a survey of some of the non-reducibility results
which apply also to $\dD2$ reducibility.

We begin with a sequence of absoluteness results which will pave the
way for forcing arguments later on.

\begin{prop}
  \label{prop_decidable_pres}
  If $A$ is an infinite time decidable set, then in any forcing
  extension there is an unambiguous interpretation of $A$ and moreover
  it remains an infinite time decidable set.
\end{prop}

\begin{proof}
  If $A$ is infinite time decidable by the program $p$ and the real
  parameter $z$, define $A$ of the forcing extension to be the set
  decided by $p$ and $z$.  To see that this is well-defined, suppose
  that programs $p,q$ both compute $A$ in the ground model. This is a
  $\dP2$ fact, and so by Shoenfield's absoluteness theorem, it remains
  true in the forcing extension.
\end{proof}

We remark that the analog of Proposition~\ref{prop_decidable_pres}
holds for infinite time eventually decidable, infinite time
semidecidable and even for absolutely $\dD2$ sets.  Similarly, we have
the following result.

\begin{prop}
  \label{prop_countable_pres}
  If an equivalence relation $E$ is infinite time enumerable, then $E$
  is infinite time enumerable in any forcing extension.
\end{prop}

\begin{proof}
  If $E$ is infinite time enumerable, then there exists an infinite
  time computable $f$ such that $f(x)$ codes $[x]_E$.  Hence we have
  that for all $x,y\in2^\omega$, the relation $x\mathrel{E}y$ holds if
  and only if there exists $n\in\omega$ such that $y=f(x)_n$.  This is
  a $\dP2$ assertion about the programs computing $f$ and deciding $E$
  and hence it is absolute to forcing extensions.
\end{proof}

\begin{prop}
  \label{prop_reduction_pres}
  Suppose that $E,F$ are absolutely $\dD2$ equivalence relations, and
  let $f$ be an infinite time eventually computable reduction from $E$
  to $F$.  Then in any forcing extension, $f$ remains such a
  reduction.
\end{prop}

\begin{proof}
  By the remarks following Proposition~\ref{prop_decidable_pres}, $E$,
  $F$, and $f$ may be unambiguously interpreted in any forcing
  extension.  Clearly, since $E,F$ are $\dD2$, the statement
  \[\forall x\forall y(x\mathrel{E}y\leftrightarrow f(x)\mathrel{F}f(y))
  \]
  is $\dP2$ and hence absolute to forcing extensions.  We must
  additionally check that $f$ remains a total function in any
  extension.  Indeed, suppose that the program $e$ eventually computes
  the function $f$.  Then $f$ is total if and only if for every
  $x\in2^\omega$ and every settled well-ordered sequence of snapshots
  according with $e$, the value $f(x)$ eventually appears on the
  output tape.  This demonstrates that the assertion ``$f$ is total''
  is $\dP2$ and hence it is absolute to forcing extensions.
\end{proof}

We remark that by Theorem~\ref{thm_semi},
Proposition~\ref{prop_reduction_pres} fails for infinite time
semicomputable reductions $f$.  On the other hand, the conclusion of
Proposition~\ref{prop_reduction_pres} does hold for absolutely $\dD2$
functions.

We now introduce several forcing methods for establishing
non-reducibility results.  Most of the arguments here have been used
in the study of Borel equivalence relations, but our adaptations will
apply even in the case of infinite time computable reductions.

\begin{defn}
  Let $E$ be any equivalence relation.  If $\PP$ is a notion of
  forcing then a $\PP$-name $\tau$ is said to be a \emph{virtual
    $E$-class} if the following hold:
  \begin{itemize}
  \item If $G$ is $\PP$-generic, then in $V[G]$,
    $\tau_G\not\mathrel{E}x$ for any $x\in V$
  \item If $G\times H$ is $\PP^2$-generic, then in $V[G\times H]$,
    $\tau_G\mathrel{E}\tau_H$
  \end{itemize}
  We say that $E$ is \emph{pinned} if it doesn't admit a virtual
  class.
\end{defn}

For instance, $E_{\textrm{set}}$ admits a virtual class via the
forcing $\PP=\Coll(\omega,\RR)$ which adds an $\omega$-sequence of
reals with finite conditions.  Any $\PP$-generic sequence will list
precisely the collection of ground model reals, and hence any two
generics will be $E_{\textrm{set}}$ equivalent.  Similarly, it is
easily seen that $\oiso_{WO}$ admits a virtual class via the forcing
$\QQ=\Coll(\omega,\omega_1)$.

The key facts, essentially due to Hjorth, are that the countable Borel
equivalence relations are pinned (see \cite[Theorem 22]{kanovei}) and
that there cannot be a Borel reduction from a non-pinned equivalence
relation to a pinned equivalence relation (see \cite[Lemma
20]{kanovei}).  Using exactly the same methods, can show the
following.

\begin{prop}
  \label{prop_pinned1}
  If $E$ is an absolutely $\dD2$ equivalence relation such that in any
  forcing extension, none of its classes can be changed by forcing,
  then $E$ is pinned.
\end{prop}

\begin{proof}
  Suppose that $E$ is not pinned, and let $\PP$ be a notion of forcing
  with a $\PP$-name $\sigma$ for a virtual class of $E$.  Let $g,h$ be
  mutually generic for $\PP$, so that $\sigma_g\mathrel{E}\sigma_h$
  holds in $V[g,h]$.  Since the classes of $E$, as interpreted in
  $V[g]$, cannot be changed by forcing, we must have $\sigma_h\in
  V[g]$.  Since $g,h$ are mutually generic, it follows that
  $V=V[g]\cap V[h]$, and hence $\sigma_h\in V$, a contradiction.
\end{proof}

For instance, any countable infinite time eventually decidable
equivalence relation satisfies the hypothesis of Proposition
\ref{prop_pinned1}.

\begin{prop}
  \label{prop_pinned2}
  Let $E$ and $F$ be absolutely $\dD2$ equivalence relations.  If
  $E\leq_eF$ and $F$ is pinned, then $E$ is pinned.
\end{prop}

\begin{proof}
  Suppose that $E$ is not pinned, and let $\PP$ be a notion of forcing
  with a $\PP$-name $\sigma$ for a virtual class of $E$.  If $f$ is an
  eventual reduction from $E$ to $F$, it is easy to see that a the
  natural $\PP$-name for $f$ applied to $\sigma$ (let us call it
  $f(\sigma)$) has the property that if $G\times H$ is
  $\PP^2$-generic, then in $V[G\times H]$, we have
  $f(\sigma)_G\mathrel{F}f(\sigma)_H$.  Since $F$ is pinned, $\PP$
  forces that $f(\sigma)\mathrel{F}y$ for some ground model real $y$.
  Now, by Shoenfield's absoluteness theorem, there exists a ground
  model real $x$ such that $f(x)\mathrel{F}y$.  It follows that $\PP$
  forces that $f(\sigma)\mathrel{F}f(x)$, and hence that $\sigma$ is
  $E$-equivalent to the ground model real $x$, a contradiction.
\end{proof}

\begin{cor}
  $E_{\textrm{set}}$ isn't eventually reducible to $E_\infty$, or even
  to $\oequiv_{e\infty}$.
\end{cor}

But recall that $\oequiv_{e\infty}$ is infinite time eventually
decidable; it is now unclear just where it should fit into the
picture.

\begin{question}
  Is $\oequiv_{e\infty}$ eventually reducible to $\oiso$?
\end{question}

We next turn to cardinality arguments.

\begin{prop}
  No infinite time computable equivalence relation which necessarily
  has $2^\omega$ many classes is eventually reducible to $\oiso_{WO}$.
\end{prop}

\begin{proof}
  Under $\neg\mathsf{CH}$, this is clear since $\oiso_{WO}$ only has
  $\omega_1$ many classes.  So just force $\neg\mathsf{CH}$ and appeal
  to Shoenfield's absoluteness theorem.
\end{proof}

\begin{cor}
  Equality $\mathord{=}$ on $2^\omega$ is not eventually reducible to
  $\oiso_{WO}$.  Hence also $E_0$, $E_\infty$, $E_{\textrm{set}}$,
  $\oiso_{HC}$, and so on, are not eventually reducible to
  $\oiso_{WO}$.
\end{cor}

The following result shows that moreover, $\oiso_{WO}$ does not
reduce to $E_\infty$.

\begin{prop}
  \label{wo_incomp_eq}
  The equivalence relation $\oiso_{WO}$ does not computably reduce to
  any equivalence relation which is necessarily countable.  In
  particular, $\oiso_{WO}$ does not reduce to any infinite time
  enumerable equivalence relation.
\end{prop}

\begin{proof}
  Let $f$ be an infinite time computable reduction from $\oiso_{WO}$
  to $E$.  Then $\im(f)$ is $\dS2$, so we may let $a$ be a parameter
  such that $S$ is $\dS2(a)$.  By the Mansfield-Solovay theorem (see
  \cite[Theorem~25.23]{jech}), if
  $\omega_1^{L[a]}<\omega_1<2^\omega$, then there is no $\dS2(a)$ set
  of size $\aleph_1$.  Since $\im(f)$ is clearly a $\dS2(a)$ set of
  size $\aleph_1$, we have reached a contradiction under these
  hypotheses.  Moreover, this situation can be forced over any model
  of \textsf{ZFC}.  Since the proposition ``the relation $\oiso_{WO}$
  doesn't computably reduce to equality'' is $\dP2$, Shoenfield's
  absoluteness theorem implies that it holds.
\end{proof}

\begin{figure}[h]
\begin{center}
\begin{tikzpicture}
\draw (0,2.5) ellipse (3 and 2.8) ;
\draw (0,3) ellipse (1.6 and 2.3) ;

\node at (0,0)  (equality) {$=$} ;
\node at (0,1)  (e0)       {$E_0$}
  edge[<-] (equality) ;
\node[circle,draw,dashed] at (-0.5,2.5) {$\oequiv_\alpha$} ;
\node[circle,draw,dashed,text width=1cm,text centered] at (6,3.5)
  {$E_\zeta$ $E_\Sigma$ $\oequiv_{e\infty}$} ;
\node at (-2,4) (ehyp)     {$E_{\textrm{hyp}}$} ;
\node at (2,3)  (jump)     {$J$} ;
\node at (0.5,4)(turing)   {$\oequiv_T$} ;
\node at (2,4)  (eqinfty)  {$\oequiv_\infty$}
  edge[<-,densely dotted] (jump) ;
\node at (0,5)  (einfty)   {$E_\infty$}
  edge[<-,densely dotted] (eqinfty)
  edge[<-] (turing)
  edge[<-] (ehyp) ;
\node at (0,6)  (eset)     {$E_{\textrm{set}}$}
  edge[<-] (einfty) ;
\node at (4,6)  (elambda)  {$E_\lambda$} ;
\node at (0,7)  (isohc)    {$\oiso_{HC}$}
  edge[<-] (eset) ;
\node at (4,7)  (isowo)    {$\oiso_{WO}$}
  edge[<-,densely dotted] (elambda);
\node at (5,7)  (eck)      {$E_{ck}$}
  edge[draw opacity=0] node {$\sim$} (isowo);
\node at (2,8)  (iso)      {$\oiso$}
  edge[<-] (isohc) edge[<-] (isowo);

\draw [->] (equality) .. controls (1.3,.3) and (2,2) .. (jump) ;
\draw [->] (equality) .. controls (2,.3) and (2.8,2.8) .. (eqinfty) ;
\end{tikzpicture}
\end{center}
\caption{Solid arrows denote computable (or better) reductions,
  dotted arrows denote eventual reductions.  The inner ellipse
  surrounds the nonsmooth countable Borel equivalence relations, the
  outer ellipse surrounds the infinite time eventually enumerable
  equivalence relations.  The dashed circles indicate open
  questions.\label{fig_diagram}}
\end{figure}
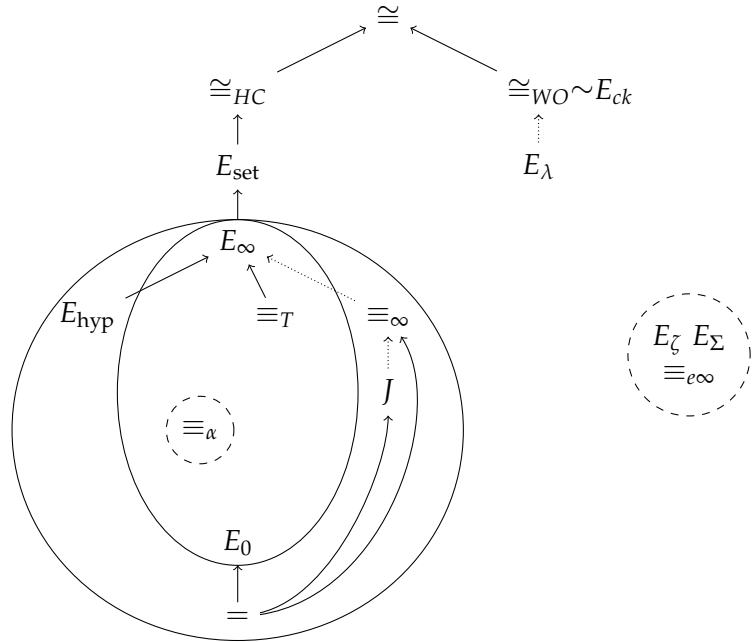

Some of the relationships between the equivalence relations considered
in this paper are summarized in Figure~\ref{fig_diagram}.

\newpage
\bibliographystyle{alphanum}
\begin{singlespace}
  \bibliography{ittm}
\end{singlespace}

\end{document}